\documentclass[11pt]{amsart}

\usepackage{amscd,amssymb,amsmath,graphicx,mathrsfs}
\usepackage[all]{xy}

\usepackage[dvips]{hyperref}
\usepackage[TS1,OT1,T1]{fontenc}

\newtheorem{theorem}{Theorem}[section]
\newtheorem{lemma}[theorem]{Lemma}
\newtheorem{corollary}[theorem]{Corollary}
\newtheorem{proposition}[theorem]{Proposition}

\theoremstyle{definition}

\newtheorem{example}[theorem]{Example}

\theoremstyle{remark}
\newtheorem{remark}[theorem]{Remark}

\numberwithin{equation}{section}

\newcommand{\C}{\mathbb C}

\newcommand{\Z}{\mathbb Z}

\newcommand{\CP}{\mathcal{P}}

\newcommand{\ca}{\mathcal}
\newcommand{\scr}{\mathscr}

\begin{document}

\title{On equivalences of derived and singular categories}
\author{Vladimir Baranovsky, Jeremy Pecharich}
\address{UC-Irvine, Mathematics Department, Irvine, CA}
\email{jpechari@math.uci.edu, vbaranov@math.uci.edu}
\label{I1}\label{I2}\label{I3}

\begin{abstract}

Let $\ca X$ and $\ca Y$ be two smooth Deligne-Mumford stacks 
and consider a pair of functions $f: \ca X \to \mathbb{A}^1$, 
$g:\ca Y \to \mathbb{A}^1$. Assuming that 
there exists a complex of sheaves 
on $\ca X \times_{\mathbb{A}^1} \ca Y$ which  
induces an equivalence of $D^b(\ca X)$ 
and $D^b(\ca Y)$, we show that there is also an equivalence of the
singular derived categories of the fibers $f^{-1}(0)$ and
$g^{-1}(0)$. 
We apply this statement in the setting of McKay correspondence, 
and generalize a theorem of Orlov on
the derived category of a Calabi-Yau hypersurface in a weighted projective
space, to products of Calabi-Yau hypersurfaces 
in simplicial toric varieties with nef anticanonical class.
\end{abstract}
\maketitle

\section{Introduction}
Mirror symmetry is a conjectural equivalence between holomorphic data on one Calabi-Yau threefold and the symplectic data on another Calabi-Yau threefold.  More precisely, Kontsevich conjectured \cite{Kontsevich:1995} given a pair of mirror Calabi-Yau threefolds $X,\widehat{X}$ there is an equivalence between the derived category of coherent sheaves on $X$ and a suitably defined Fukaya category on $\widehat{X}$.  Later this was extended to involve Fano varieties.  The mirror to a Fano variety, $X$ is no longer just a variety but a pair $(\widehat{X},W)$ where $\widehat{X}$ is a regular scheme and $W:\widehat{X}\rightarrow \mathbb{A}^1$ is a regular  morphism.  Such a pair is called a Landau-Ginzburg model in the physics literature.  When $\widehat{X}$ is affine Kontsevich suggested replacing the derived category of coherent sheaves with a category of $2$-periodic complexes where the composition is no longer $0$ but multiplication by $W$.  Orlov defined this category and showed it is a triangulated category \cite{Orlov:2004}.  This category is denoted by $DB(W)$ and mirror symmetry is now conjectured to be
an equivalence between $DB(W)$ and a Fukaya type category.

The bounded derived category of coherent sheaves $D^b(X)$ has a triangulated subcategory $\mathfrak{Perf}(X)$ consisting of perfect complexes i.e. complexes which are locally quasi-isomorphic to a bounded complex of locally free sheaves.     If $X$ is non-singular then every bounded complex of coherent sheaves 
admits a locally free resolution.  This means that $\mathfrak{Perf}(X)$ is equivalent to $D^b(X)$.  When $X$ is singular this is no longer true.  Orlov \cite{Orlov:2004} introduced the singular derived category $D_{sg}(X)$ as the quotient of $D^b(X)$ by the full triangulated subcategory $\mathfrak{Perf}(X)$. He showed
that for affine $X$ the category $DB(W)$ is equivalent to the 
product of $D_{sg}(X_w)$ over the critical values $w$ of $W: X \to \mathbb{A}^1$. Thus, one can say that $DB(W)$ reflects the singularities of the 
fibers of $W$.

In this paper we study $D^b(\cdot)$ and $D_{sg}(\cdot)$ for 
Deligne-Mumford stacks. 

\bigskip
\noindent
\textit{Convention}.  The term \textit{stack}  will mean a
Deligne-Mumford stack $\ca X$ of finite type over a field $k$ of 
characteristic zero ($k = \mathbb{C}$ in applications)
which has finite stabilizers. In addition we will always assume that
$\ca X$ is \textit{quasi-projective}, i.e. admits a locally closed 
embedding
in a smooth Deligne-Mumford stack $\ca W$, such that 
$\ca W$ is proper over $Spec\;k$
and has a projective coarse moduli space. 
If we can choose a closed embedding into such $\ca W$ 
then $\ca X$ is called \textit{projective}.
Quasi-projectivity ensures that $\ca X$ 
has a quasi-projective moduli space and that every coherent sheaf on $\ca X$
is a quotient of a vector bundle (the so-called \textit{resolution property}).
 See \cite{Kr}
for an excellent overview of related results. In particular, if $\ca X$
is smooth then every coherent sheaf admits a finite locally free resolution. 
For a quasi-projective stack $\ca X$, we denote by $D^b(\ca X)$ the 
bounded derived
category of coherent sheaves. All functors between such 
categories (pullbacks, pushforwards and tensor products) are assumed to 
be derived (otherwise they  are not well defined on $D^b$) hence we drop the usual
letters $R$ and $L$ from the notation.

\bigskip
\noindent
The main technical result of this paper, proved in Section 2, is as follows.

\begin{theorem} 
Let $\ca X$, $\ca Y$ be two smooth  stacks, 
$f: \ca X \to \mathbb{A}^1$, $g: \ca Y\to \mathbb{A}^1$ two morphisms
and $\widetilde{\scr F}$ a complex on $\ca X \times_{\mathbb{A}^1}
\ca Y$ with support proper over $\ca X$ and $\ca Y$.
Let $\Phi_{\scr F}: D^b(\ca X) \to D^b(\ca Y)$ be the Fourier-Mukai transform
with the kernel $\scr F = h_* \widetilde{\scr F}$, where $h:
\ca X \times_{\mathbb{A}^1} \ca Y \to \ca X \times \ca Y$ is the closed
embedding. If $\Phi_{\scr F}$ is an equivalence of categories, 
and $\ca X_0$, $\ca Y_0$ are the fibers of $f$, $g$ over $0 \in \mathbb{A}^1$, 
respectively, then the pullback of $\widetilde{\scr F}$ 
to $\ca X_0 \times \ca Y_0$ induces an
equivalence of derived categories $D^b(\ca X_0) \simeq D^b(\ca Y_0)$ 
which descends to an equivalence of singular derived categories
$D_{sg}(\ca X_0) \to D_{sg}(\ca Y_0)$.
\end{theorem} 

\bigskip
\noindent
This theorem can be applied to the situation when one has 
the following diagram of  stacks and proper birational morphisms 
 
\begin{equation}
\label{main}
\xymatrix{                   &\ca W\ar[ld]_\mu \ar[rd]^\nu\\
                  \ca X\ar[rd]_{\phi}&  &\ca Y\ar[ld]^{\psi}\\
                                   &Z}
\end{equation}  
We assume that  $\ca X, \ca Y, \ca W$ are smooth
stacks, $Z$
is an affine variety, 
$\mu^*(K_{\ca X})=\nu^*(K_{\ca Y})$ and in addition $\phi$ and 
$\psi$ induce projective morphisms of coarse moduli spaces. 
For instance, we can assume that $Z$ is a  singular quotient
$V/\Gamma$ of a finite dimensional vector space $V$ by a linear action of
a finite group $\Gamma \subset SL(V)$, 
$\ca X = [V/\Gamma]$ is the same quotient
considered as a smooth stack and $\ca Y$ is a crepant 
resolution of the variety $V/\Gamma$.

By a strengthened version of derived McKay Correspondence conjecture,
such a diagram \eqref{main} should imply existence of a Fourier-Mukai
equivalence with a kernel obtained by direct image from
$\ca X \times_Z \ca Y$. In particular, $\widetilde{\scr F}$ exists
whenever the morphisms $f: \ca X \to \mathbb{A}^1$, $g: \ca Y \to \mathbb{A}^1$
are obtained by pulling back the same regular function $Z \to \mathbb{A}^1$. 
This holds  in the cases of the derived McKay correspondence considered
by Bezrukavnikov and Kaledin; Bridgeland, King, and Reid; and Kawamata \cite{BezKaledin:2004,BKR:2001,Kawamata:2005}.
Thus, our first  application, Theorem 2.10, is in these three cases.

Our Theorem 2.10 generalizes an earlier result of  
Quintero-V\'elez, cf. \cite{Velez:2008}, who proves the statement under the
assumption that $\ca Y$ is given by $G-Hilb(V)$, the Hilbert scheme
of $G$-clusters, and that $\Phi$ is given by the structure sheaf of
the universal subscheme. In particular, case (1) of Theorem 2.10
follows from \textit{loc. cit.}. As observed in 
the same paper, the assumption also  holds when 
$\Gamma = \mathbb{Z}/n\mathbb{Z}$ where $n = \dim V$,
$\Gamma$ acts diagonally on $V$ and $\ca Y$ is the total
space of the canonical bundle on the projective space $\mathbb{P}(V)$. 
After the first version of this paper appeared as preprint, we 
have learned about an even earlier work by  Mehrotra, 
cf. \cite{Me}, who have constructed a full and faithful embedding 
$D_{sg}(\ca X_0) \to D_{sg}(\ca Y_0)$ under the same assumption
(i.e. that $\ca X = [V/\Gamma]$ and $\ca Y = G-Hilb(V)$). 
We also mention here that in the case of schemes, 
but for arbitrary base change and without any flatness
assumption, the results of Section 2.2 
were proved earlier by Kuznetsov, cf. Section 2.7 of \cite{Ku}.

In Section 3 we give another application in the setting when 
all four spaces and morphisms are toric with respect to the action
of the same (split) torus $T$. 
We apply Theorem 1.1 and results of Kawamata, cf. \cite{Kawamata:2005}, 
to the case when $\Gamma \subset SL(V)$
is abelian and $f$ is given by a $\Gamma$-invariant polynomial 
on $V$.   
We mention here that the work \cite{HLdS} on Fourier-Mukai transform 
for Gorenstein schemes also deals with a similar situation, although 
from a somewhat different angle. To formulate the second application,
cf. Corollaries 3.3 and 3.7, first 
let $\mathbb{P}(\overline{a}):=\mathbb{P}(a_0,\cdots,a_n)$ with $a_i>0$ for all $i$, be the weighted projective space.  Given $f$ a quasi-homogeneous polynomial that is invariant under the action of $\Z_N$ where $N=\sum a_i$ the zero set $Y$ is a Calabi-Yau hypersurface of $\mathbb{P}(\overline{a})$.  Orlov gave an algebraic proof of an equivalence $D^b(Y)\cong D_{sg}^{\C^*}(\C^{n+1},f)$ \cite[Thm 3.12]{Orlov:2005}, where $\C^*$ acts on $\C^{n+1}$ with weights $(a_0, \ldots, 
a_n)$ and $D_{sg}^{\C^*}$ stands for the equivariant version of the
singular category.
As a corollary of Theorem 1.1 and Kawamata's theorem on toric crepant resolutions \cite[Prop. 4.2]{Kawamata:2005} we give a geometric proof of this statement and also generalize it to products of Calabi-Yau 
hypersurfaces in simplicial toric varieties with nef anticanonical class.

\bigskip
\noindent
\textit{Acknowledgements}.  We would like to thank Tony Pantev for valuable comments. The work of the first author was partially supported by 
the Sloan Research Fellowship.

\section{Proof}

In this section we consider the following commutative diagram in which
all horizontal arrows are regular closed embeddings of codimension one, and
all other arrows are the canonical projections:
\begin{equation}
\label{diagram}
\xymatrix{
\ca Y_0 \ar[r]^{i_0} & \ca Y \\
\ca X_0\times \ca Y_0\ar[r]^{k_0}\ar[d]_{\pi_{\ca X_0}} \ar[u]^{\pi_{\ca Y_0}}
&\ca X\times_{\mathbb{A}^1}\ca Y\ar[d]_{p_{\ca X}}\ar[r]^h \ar[u]^{p_{\ca Y}}
&\ca X\times \ca Y\ar[ld]^{\pi_{\ca X}} \ar[lu]_{\pi_{\ca X}}\\
\ca X_0\ar[r]_{j_0}&\ca X}
\end{equation}
Note that the both squares are cartesian and the vertical arrows represent
flat morphisms.  We will also consider the shifted line bundles
\begin{equation}
\label{om}
\omega := \pi^*_{\ca X}K_{\ca X} [n], \qquad
\widetilde{\omega}:= h^*\omega[-1], \qquad
\omega_0 = \C^*_0 \widetilde{\omega}
\end{equation}
on $\ca X \times \ca Y$, $\ca X \times_{\mathbb{A}^1}\ca Y$ and 
$\ca X_0 \times \ca Y_0$, respectively.

\subsection{Generalities on sheaves and stacks.}

To prove Theorem 1.1 we first need a few lemmas. The reader may wish to
skip to Section 2.2 and return to the statements of this section 
as they are quoted in the proof.

\begin{lemma} \cite[Prop. 13.1.9]{Champs} Let 
\begin{equation*}
\xymatrix{                   \ca Y' \ar[d]_\xi \ar[r]^v&\ca Y\ar[d]^\zeta\\
                  \ca X'\ar[r]_u&  \ca X}\end{equation*}
                  be a cartesian diagram of Deligne-Mumford stacks with $\zeta$ quasi-compact and $u$ flat and $\scr N\in D^+_{qc}(\ca Y)$. Then the canonical morphism 
                  \begin{equation*} u^*\zeta_*\scr N\rightarrow \xi_*v^*\scr N \end{equation*} 
is an isomorphism in $D_{qc}^+(\ca X')$.  
                
  \end{lemma}

\begin{lemma} Let $\ca X$ be a smooth quasi-projective stack. There exists
a smooth projective stack $\overline{\ca X}$ and an open 
embedding $\ca X \to \overline{\ca X}$. 
\end{lemma}
\begin{proof}
By  assumption $\ca X$ admits a locally closed embedding
in a smooth  projective stack $\ca W$. By \cite{Kr} we 
can assume that $\ca W = P/GL(n)$ where $P$ is a smooth quasi-projective 
variety and $GL(n)$ acts on $P$ linearly and with 
finite stabilizers. Then also $\ca X \simeq Q/GL(n)$ where $Q\subset P$ is
smooth and locally closed. The closure $\overline{Q}$ of $Q$ in $P$ is
$GL(n)$-invariant, although it may be singular. Using canonical desingularization (i.e. versions due to
Bierston-Millman, Villamayor and W\l{}odarczyk) by 
blowing up a maximal stratum of a certain local invariant which is automatically 
$GL(n)$-invariant, we can 
find an iterated $GL(n)$-equivariant blowup $\widetilde{P} \to P$ with smooth 
$GL(n)$-invariant closed centers, such that the proper transform $\widetilde{Q}$
is smooth and the morphism $\widetilde{Q} \to \overline{Q}$ 
 restricts to an isomorphism over $Q$. Then
$\widetilde{Q}$ is quasi-projective, and the stabilizers of the $G$-action
on $\widetilde{Q}$ are still finite (since they embed into the stabilizers
of points in $\overline{Q}$). Moreover, the moduli space $\widetilde{Q}/GL(n)$
 of the quotient stack
$[\widetilde{Q}/GL(n)] =:\overline{\ca X}$ is projective. In fact, 
since $\widetilde{Q}/GL(n)$ is closed in
$\widetilde{P}/GL(n)$, it suffices to show that the morphism 
$\widetilde{P}/GL(n) \to P/GL(n)$ is projective. By induction we 
can assume that $\widetilde{P}$ is a single blowup of $P$ at a smooth
$GL(n)$-invariant center $R$. We can find an ample line bundle 
$L$ on $P$ and a finite-dimensional $G$-invariant subspace
of sections $U \subset \Gamma(P, L)$ such that the common zero scheme
of these sections is  $R$. Then $\widetilde{P}$ can 
be identified with a closed subvariety of $P \times \mathbb{P}(V^*)$ 
(i.e. closure of the graph of the rational map defined by the linear system
$|U|$)
thus it suffices to show that $P \times \mathbb{P}(V^*)/GL(n)
\to P/GL(n)$ is a projective morphism, which follows easily by a
GIT-type argument.
\end{proof}
The next lemma has a relatively quick proof due to the  
(quasi)-projectivity condition which imposed on  stacks.
\begin{lemma} (Relative Serre Duality: smooth projective case) 
Let $\ca X$ be a smooth projective stack of dimension $n$ and $\ca Y$
a quasi-projective stack. The functor $\pi_{\ca Y*}:
D^b(\ca X \times \ca Y) \to D^b(\ca Y)$ has a right adjoint 
$$
\pi^!_{\ca Y}(\cdot) \simeq \pi^*_{\ca X}K_{\ca X}[n] \otimes \pi^*_{\ca Y} (\cdot)
$$
\end{lemma}
\begin{proof}
We first observe that
our stacks  satisfy the resolution 
property and morphisms are separated. Hence by Proposition 1.9 in 
\cite{Ni} the right adjoint  $\pi^!_{\ca Y}$ of $\pi_{\ca Y*}$
exists although apriori it is defined as a functor 
$D(\ca Y) \to D(\ca X \times \ca Y)$
on unbounded derived categories of complexes of $\mathcal{O}$-modules with
quasi-coherent cohomology. By base change Lemma 2.1, existence of 
$\pi^!_{\ca Y}$ and the proof of Theorem 5.4 in \cite{Ne} we can conclude 
that 
$$
\pi^!_{\ca Y}(\cdot) \simeq 
\pi^!_{\ca Y}(\mathcal{O}_{\ca Y}) \otimes \pi^*_{\ca Y} (\cdot).
$$
It remains to establish 
$$
\pi^!_{\ca Y}(\mathcal{O}_{\ca Y}) \simeq \pi_{\ca X}^* K_{\ca X}[n]
$$
When $\ca Y \simeq Spec(k)$, this is the contents of
Theorem 1.32 (Smooth Serre Duality) in \cite{Ni}. For general 
$\ca X$ we consider the diagram 
\begin{equation*}\xymatrix{
\ca X\times \ca Y\ar[r]^{\pi_{\ca X}}\ar[d]_{\pi_{\ca Y}}&\ca X\ar[d]^{p}\\
\ca Y\ar[r]_{q}& Spec(k)}
\end{equation*}
Since $p^!\mathcal{O}_{Spec(k)} \simeq K_{\ca X}[n]$ by \textit{loc. cit.}
it remains to show that $\pi_{\ca X}^* p^! \simeq \pi^!_{\ca Y} q^*$.
Following Verdier, we consider the composition 
$$
\pi_{\ca X}^* p^! \to \pi^!_{\ca Y}\pi_{\ca Y*}\pi_{\ca X}^* p^! 
\simeq \pi^!_{\ca Y} q^* p_* p^! \to \pi^!_{\ca Y} q^*,
$$
where we have used the morphisms induced by adjunctions
 $Id \to \pi^!_{\ca Y}\pi_{\ca Y*}$, 
$p_* p^! \to Id$ and the base change formula from Lemma 2.1.
The morphism of functors $\pi_{\ca X}^* p^! \to \pi^!_{\ca Y} q^*$
is functorial with respect to $\ca Y$ hence in proving that it is an isomorphism
it suffices to replace $\ca Y$ by a scheme $Z$ admitting a
finite flat surjective morphism $Z \to \ca Y$, which exists by \cite{Kr}.
Now we apply Theorem 1.23 in \cite{Ni}. This finishes the proof.
\end{proof}

\begin{corollary}
Recall the notation of \eqref{diagram} and \eqref{om} and let $\scr H$, 
resp. $\scr H_0$ be a complex on $\ca X \times \ca Y$, resp.  
$\ca X_0 \times \ca Y_0$, which has proper support over $\ca Y$, resp. 
$\ca Y_0$. Let also $b$, resp. $b_0$, be a complex on $\ca Y$, resp. $\ca Y_0$. Then there are isomorphisms of bifunctors
\begin{equation}
\label{supp-dual}
Hom_{\ca Y} (\pi_{\ca Y*} \scr H, b) \simeq Hom_{\ca X \times \ca Y}
(\scr H, \omega \otimes \pi^*_{\ca Y}(b))
\end{equation}
\begin{equation}
\label{supp-dual-zero}
Hom_{\ca Y_0} (\pi_{\ca Y_0*} \scr H_0, b_0) \simeq Hom_{\ca X_0 \times \ca Y_0}
(\scr H, \omega_0 \otimes \pi^*_{\ca Y_0}(b_0)).
\end{equation}
\end{corollary}
\begin{proof}
For the first isomorphism, choose a smooth compactification $\ca X \to \overline{\ca X}$
as in Lemma 2.2, view $\scr H$ as an object on $\overline{\ca X} \times \ca Y$
and then apply Lemma 2.3. 

For the second isomorphism consider the upper cartesian 
square in \eqref{diagram}. Replacing $\ca X$ by a larger open substack in the
compactification $\overline{\ca X}$
(and possibly blowing up $\overline{\ca X}$ but without changing $\ca X$) 
we can assume that $f: \ca X \to 
\mathbb{A}^1$ is proper (then $\scr H_0$ will be replaced by its direct 
image to this larger open substack). Therefore 
 we can assume that $p_{\ca Y}$ and $\pi_{\ca Y_0}$ are
proper. As in the end of the proof 
of Lemma 2.3 we have $k_0^* p_{\ca Y}^! \simeq
\pi_{\ca Y_0}^! i_0^*$. Then 
$$
Hom_{\ca Y_0} (\pi_{\ca Y_0*} \scr H_0, b_0) \simeq Hom_{\ca X_0 \times \ca Y_0}
(\scr H_0, \pi_{\ca Y_0}^!(b_0))
\simeq Hom_{\ca X_0 \times \ca Y_0}
(\scr H_0, \pi^!_{\ca Y_0} (\mathcal{O}_{\ca Y_0}) \otimes 
\pi^*_{\ca Y}(b_0))
$$
and we only need to apply 
$$
\pi_{\ca Y_0}^! (\mathcal{O}_{\ca Y_0})
\simeq \pi_{\ca Y_0}^! i_0^* (\mathcal{O}_{\ca Y})
\simeq k_0^* p_{\ca Y}^!(\mathcal{O}_{\ca Y}) \simeq \omega_0.
$$
\end{proof}
\begin{lemma}
The complex $\scr F_0 = k_0^* \widetilde{\scr F}$ is bounded and has proper support over $\ca X_0$ and $\ca Y_0$.
Moreover, the following properties hold for 
$\scr F_0$:
\begin{enumerate}
\item for any perfect complex
$\scr K$ on $\ca X_0 \times \ca Y_0$ the direct image
$\pi_{\ca X_0*}(\scr F_0 \otimes \scr K)$, resp. 
$\pi_{\ca Y_0*}(\scr F_0 \otimes \scr K)$ is a perfect complex on 
$\ca X_0$, resp. $\ca Y_0$;

\item for any complex $\scr H$ in
$D^b(\ca X_0)$, resp. $D^b(\ca Y_0)$, the complex 
$\scr F_0 \otimes \pi^*_{\ca X_0} \scr H$, resp. $\scr F_0 \otimes \pi^*_{\ca Y_0}
\scr H$ is an object of $D^b(\ca X_0 \times Y_0)$.
\end{enumerate}
\end{lemma}
\begin{proof}
The support property is obvious since we assume that $\widetilde{\scr F}$
has proper support over $\ca X$ and $\ca Y$. Boundedness of $\scr F_0$ 
holds since $k_0$ is a regular closed embedding of codimension one, hence
has finite tor dimension. As for the properties (1) and (2), it suffices
to prove those which involve $\pi_{\ca X_0}$.

Assume that $\widetilde{\scr F}$ has cohomology in degrees $[-n, 0]$.
Since $\ca X \times_{\mathbb{A}^1} \ca Y$ has the resolution property,
we can construct a complex of vector bundles on it
$$
E_{-N} \to E_{-N+1} \to \ldots \to E_{-1} \to E_0 \to 0
$$
with $N > n + \dim \ca X$   
which is quasi-isomorphic to $\widetilde{\scr F}$ except possibly in 
degree $-N$. By a standard argument using the smoothness of
$\ca X$, the kernel of $E_{-N} \to E_{-N+1}$
is projective over $\mathcal{O}_{\ca X}$. Thus $\widetilde{\scr F}$
is quasi-isomorphic to a finite complex of coherent sheaves which 
are projective over $\ca X$. Similarly, its pullback $\scr F_0$ is
isomorphic to a finite complex of coherent sheaves which are
projective over
$\mathcal{O}_{\ca X_0}$. This proves (2).

To show that $\pi_{\ca X_0*}(\scr F_0 \otimes \scr K)$ is perfect is 
suffices to show that $\pi_{\ca X_0*}(\scr F_0 \otimes \scr K) \otimes 
\scr H$ is bounded for any bounded $\scr H$ in $D^b(\ca X_0)$. 
By projection formula it suffices to show that
$\scr F_0 \otimes \scr K \otimes \pi_{\ca X_0}^* \scr H$
is bounded which is immediate from the finite  
$\mathcal{O}_{\ca X_0}$-projective resolution of $\scr F_0$.
\end{proof}

\subsection{Main argument}

\bigskip
\noindent
Lemma 2.5 implies that we have a well-defined Fourier-Mukai 
transform 
$F_0:=\pi_{\ca Y_0*}(\scr F_0\otimes \pi_{\ca X_0}^*(\cdot)):D^b(\ca X_0)\rightarrow D^b(\ca Y_0)$, cf. \cite{Huy}. There is also a similar Fourier-Mukai
transform $F: D^b(\ca X)\to D^b(\ca Y)$ with the kernel 
$\scr F = h_*\widetilde{\scr F}$ (note that $\scr F$ is a perfect complex
since $\ca X$ and $\ca Y$ are smooth). 

\begin{lemma} Let $i_0:\ca Y_0\rightarrow \ca Y$ and 
$j_0:\ca X_0\rightarrow \ca X$ be the closed immersions of the fibers then 
there is a functorial isomorphism
\begin{equation*} i_{0*}F_0\cong Fj_{0*} \end{equation*} \end{lemma}

\begin{proof}
Similarly to Theorem 6.1 in \cite{Ch} we use a series of isomorphisms 

\begin{align*}  
Fj_{0*}(\cdot)&=\pi_{\ca Y*}(\pi^*_{\ca X}(j_{0*}(\cdot))\otimes \scr G)\\
&=\pi_{\ca Y*}(\pi_{\ca X}^{*}(j_{0*}(\cdot))\otimes j_*\scr F)\\
&=\pi_{\ca Y*}(h_*(h^*\pi_{\ca X}^*j_{0*}(\cdot)\otimes \scr F))\\
&=\pi_{\ca Y*}(h_*(p^*_{\ca X}j_{0*}(\cdot)\otimes \scr F))\\
&=\pi_{\ca Y*}(h_{*}(k_{0*}\pi_{\ca X_0}^*(\cdot)\otimes \scr F))\\
&=\pi_{\ca Y*}(h_*k_{0*}(\pi_{\ca X_0}^*(\cdot)\otimes k_0^*\scr F))\\
&=\pi_{\ca Y*}(i_0\times j_0)_*(\pi_{\ca X_0}^*(\cdot)\otimes \scr F_0)\\
&=i_{0*}\pi_{\ca Y_0*}(\pi_{\ca X_0}^*(\cdot)\otimes \scr F_0)\\
&=i_{0*}F_0(\cdot)
\end{align*}
The third and sixth isomorphisms are due to the projection formula, forth is by commutativity of the lower triangle in \eqref{diagram}, 
and the fifth is the base change isomorphism of Lemma 2.1 
applied to the lower square of the same diagram.
\end{proof}

\begin{lemma}
In the notation of \eqref{diagram}, \eqref{om},
define a complex in $D^b(\ca X \times_{\mathbb{A}^1} \ca Y)$:
$$
\widetilde{\scr G} = \mathcal{H}om(\widetilde{\scr F},  \widetilde{\omega}).
$$
Then
the right adjoint $G: D^b(\ca Y) \to D^b(\ca X)$ to $F$ is given by 
the Fourier-Mukai transform 
with the kernel $\scr G = h_* \widetilde{\scr G}$.
Similarly, the right adjoint
$G_0: D^b(\ca Y_0) \to D^b(\ca X_0)$ to $F_0$ is given by the 
Fourier-Mukai transform with the kernel $\scr G_0 = k_0^* \widetilde{\scr G}$.
\end{lemma}
\begin{proof}
We first note that boundedness of $\widetilde{\scr G}$ follows by 
repeating the argument of Lemma 2.5 and using smoothness of $\ca X$.

To prove the assertion about $F$, assume that $a$, resp. 
$b$, is a object in $D^b(\ca X)$, resp. $D^b(\ca Y)$, and
recall that we denoted $\scr F = h_* \widetilde{\scr F}$. 
 Therefore by \eqref{supp-dual}
\begin{align*} 
 Hom_{\ca Y} (F(a), b) 
& =  Hom_{\ca Y}(\pi_{\ca Y*}(\scr F \otimes \pi_{\ca X}^*(a)), b) 
& \simeq  Hom_{\ca X \times  \ca Y} (\scr F \otimes \pi^*_{\ca X}(a), 
\omega \otimes \pi^*_{\ca Y}(b))\\
&\simeq  Hom_{\ca X \times  \ca Y} (\pi^*_{\ca X}(a), 
\scr F^\vee \otimes \omega\otimes \pi^*_{\ca Y}(b))
&\simeq  Hom_{\ca X} (a,\pi_{\ca X*} (\scr F^\vee \otimes \omega
\otimes \pi^*_{\ca Y}(b)))
\end{align*}
In the second line we use the fact that $\scr F$ is perfect
and  the adjunction between 
$\pi^*_{\ca X}$ and $\pi_{\ca X*}$.
Observe that the second argument of the last expression is a 
Fourier-Mukai transform of $b$ with the kernel 
$\scr F^\vee \otimes \omega \simeq \mathcal{H}om(\scr F, \omega)$. Thus
by definition on $\scr F$ and $\widetilde{\omega}$
it suffices to show that
$$
\mathcal{H}om(h_* \widetilde{\scr F}, \omega)
\simeq h_*\mathcal{H}om(\widetilde{\scr F}, h^* \omega[-1])
$$
Since $h^!(\cdot) \simeq h^*(\cdot)[-1]$ this follows from 
Corollary 1.22 of \cite{Ni} (we may even reduce to the case
of schemes constructing a morphism from the LHS to the RHS
as in \textit{loc. cit.} and then checking that it is an isomorphism on
etale local affine charts). This finishes the proof for the functor $G$.

\bigskip
\noindent
The assertion about $G_0$ is proved similarly: let $\scr H_0 
= \scr F_0 \otimes \pi^*_{\ca X_0}(b)$ and apply 
\eqref{supp-dual-zero} to obtain
\begin{align*} 
Hom_{\ca Y_0} (F_0(a), b) &
= Hom_{\ca Y}(\pi_{\ca Y_0*}(\scr F_0 \otimes \pi_{\ca X_0}^*(a)), b) \\
& \simeq Hom_{\ca X_0 \times \ca Y_0} (\scr F_0 \otimes \pi_{\ca X_0}^*(a), 
\omega_0 \otimes \pi_{\ca Y_0}^* (b)) \\
& \simeq Hom_{\ca X_0 \times \ca Y_0} (\pi_{\ca X_0}^*(a), 
\mathcal{H}om(\scr F_0, \omega_0 \otimes \pi_{\ca Y_0}^* (b)))
\end{align*}
Repeating the proof of Lemma 3.5 in \cite{Ba} we derive
from our Lemma 2.5 that
$$
\mathcal{H}om(\scr F_0, \omega_0 \otimes \pi_{\ca Y_0}^* (b))
\simeq \mathcal{H}om(\scr F_0, \omega_0) \otimes \pi_{\ca Y_0}^* (b),
$$
thus our assertion reduces to 
$$
k_0^* \mathcal{H}om(\widetilde{\scr F}, \widetilde{\omega}) \simeq 
\mathcal{H}om(\scr F_0, \omega_0) = \mathcal{H}om(\C^*_0\widetilde{\scr F}, 
k_0^* \widetilde{\omega}).
$$
Since $\widetilde{\omega}$ is a shift of a line bundle, the last isomorphism
follows immediately by replacing $\widetilde{\scr F}$ with 
an $\mathcal{O}_{\ca X}$-projective resolution as in the proof of
 Lemma 2.5.
\end{proof}

\begin{lemma}
With the notation as in Lemma 2.7 there is a functorial isomorphism
\begin{equation*} j_{0*}G_0\cong Gi_{0*} \end{equation*}
\end{lemma}

\begin{proof} The proof is exactly the same as Lemma 2.6 therefore we omit it.  \end{proof}

To finish the proof of Theorem 1.1 we need a category theory lemma:
\begin{lemma}\cite[Lemma 1.2]{Orlov:2004} Let $\mathscr{N}$ and $\mathscr{N}'$ be full triangulated subcategories of triangulated categories $\mathscr{D}$ and $\mathscr{D}'$ respectively.  Let $F:\mathscr{D}\rightarrow \mathscr{D}'$ and $G:\mathscr{D}'\rightarrow \mathscr{D}$ be an adjoint pair of exact functors such that $F(\mathscr{N})\subset \mathscr{N}'$ and $G(\mathscr{N}')\subset \mathscr{N}$.  Then they induce functors 
\begin{equation*} \overline{F}:\mathscr{D}/\mathscr{N}\rightarrow \mathscr{D}'/\mathscr{N}'\hspace{0.2in} \overline{G}:\mathscr{D}'/\mathscr{N}'\rightarrow \mathscr{D}/\mathscr{N} \end{equation*}
which are adjoints.  Moreover, if the functor $F:\mathscr{D}\rightarrow \mathscr{D}'$ is fully faithful, the functor $\overline{F}:\mathscr{D}/\mathscr{N}\rightarrow \mathscr{D}'/\mathscr{N}'$ is also fully faithful.  \end{lemma}

\noindent \textit{Proof of Theorem 1.1}.  Let $\scr E\in D^b(\ca X_0)$ and consider the exact triangle
\begin{equation*}  \scr E\rightarrow G_0F_0\scr E\rightarrow \scr C \rightarrow \scr E[1] \end{equation*}
where $\scr C$ is a cone over the morphism $\scr E\rightarrow G_0F_0\scr E$.  Applying $j_{0*}$ yields an exact triangle
\begin{equation*} j_{0*}\scr E\rightarrow j_{0*}G_0F_0\scr E\rightarrow j_{0*}\scr C\rightarrow j_{0*}\scr E[1] \end{equation*}
 But $j_{0*}G_0F_0\cong Gi_{0*}F_0\cong GFj_{0*}$
where the first isomorphism is by Lemma 2.8 and the second by Lemma 2.6.  
Since $GF$ is isomorphic to the identity morphism and $j_{0}$ is a closed immersion, $\ca C\cong 0$.  This implies that $\scr E\rightarrow G_0F_0\scr E$ is an isomorphism.

By Lemma 2.5 the functor $F_0:=\pi_{\ca Y_0}(\scr F_0\otimes \pi^*_{\ca X_0}(\cdot))$ takes perfect complexes to perfect complexes, and
similarly for $G_0$.   By Lemma 2.7  these
functors induce a pair of adjoint functors $\overline{F}_0:D_{sg}(\ca X_{0})\rightarrow D_{sg}(\ca Y_0)$ and $\overline{G}_0:D_{sg}(\ca Y_0)\rightarrow D_{sg}(\ca X_0)$. Moreover, the composition $\overline{G}_0\overline{F}_0$ is isomorphic to the identity on $D_{sg}(\ca X_0)$ and similarly 
$\overline{F}_0\overline{G}_0$ is isomorphic to the identity on $D_{sg}(\ca Y_0)$. \qed

\subsection{Applications to McKay correspondence}

In this subsection we assume that $V$ is a finite dimensional vector space
of $k = \mathbb{C}$ and $\ca X$ is the quotient stack $[V/\Gamma]$, where
$\Gamma \subset SL(V)$ is a finite subgroup acting
freely in codimension one. We take $\ca Y$ to be
a crepant resolution of the singular quotient variety $Z= V/\Gamma$ (thus,
we implicitly assume that $\ca Y$ exists which is not always the case). 
Finally, let $t: Z \to \mathbb{A}^1$ be a morphism inducing the compositions
$f: \ca X\to \mathbb{A}^1$, $g: \ca Y\to \mathbb{A}^1$ and denote by 
$\ca X_0$, $\ca Y_0$ the fibers of $f$ and $g$ over $0 \in \mathbb{A}^1$, 
respectively.

\begin{theorem}
In the above setting, suppose that there exists a complex
$\overline{\scr F}$ in $D^b(\ca X \times_Z \ca Y)$ such that
its direct image $\scr F$ to $\ca X \times \ca Y$ induced a Fourier-Mukai
transform which is an equivalence. Let $\scr F_0 = k_0^* j_* \overline{\scr F}$
where $j: \ca X \times_Z \ca Y \to \ca X \times_{\mathbb{A}^1} \ca Y$
and $k_0: \ca X_0 \times \ca Y_0 \to \ca X \times_{\mathbb{A}^1} \ca Y$
are the natural closed embeddings. Then the Fourier-Mukai transform with 
the kernel $\scr F_0$ induces equivalences
$$
D^b(\ca X_0) \simeq D^b(\ca Y_0), \quad 
\mathfrak{Perf}(\ca X_0) \simeq \mathfrak{Perf}(\ca Y_0), 
\quad D_{sg} (\ca X_0) \simeq D_{sg}(\ca Y_0)
$$
In particular, such equivalences hold either of the three cases:
\begin{enumerate}
\item $\dim V = 2$ or $3$ and $\ca Y = G-Hilb$
\item $V$ is a symplectic vector space and $\Gamma \subset Sp(V)$
\item $\Gamma$ is a finite abelian subgroup of $SL(V)$  and 
$\ca Y$ is projective over $Z$.
\end{enumerate}
\end{theorem}
\begin{proof}
The first part of the statement is an immediate consequence of Theorem 1.1.
In the second part, we only need to establish existence of appropriate $\widetilde{\scr F}$. For (1) and (2) this follows from the main results in 
\cite{BKR:2001} and \cite{BezKaledin:2004}, respectively.
For (3) we observe that since $\Gamma$ is abelian, 
there exists an $n$-dimensional torus $T$ such that $Z$ is a toric variety
with the torus $T$, and that by the proof of Corollary 3.5 in
\cite{K3} the crepant resolution $\ca Y$ is automatically toric with the
same torus $T$. Note that the corollary quoted is stated for 
projective varieties, and in our case we need to repeat its proof
invoking the relative toric version of MMP from 
\cite{FS}. Existence of $\widetilde{\scr F}$ is proved in 
Theorem 3.1 by a simple extension of results in 
\cite{Kawamata:2005}.
\end{proof}
\begin{remark}
We could require that the action of $\Gamma$ is
generically free (i.e. stabilizers \textit{are} allowed in codimension one). 
In this case consider the normal
subgroup $\Pi \subset \Gamma$ generated by complex reflections
in $V$. Then the quotient
$V' = V/\Pi$ is a smooth variety and in Theorem 2.10 one can
consider $\ca X = [V'/\Gamma']$ with $\Gamma' = \Gamma/\Pi$.
Alternatively, we can keep $\ca X = V/\Gamma$ but view $\ca Y$ not
as a smooth variety but as a smooth stack with non-trivial cyclic stabilizers
in codimension one (see the beginning of the next section).
\end{remark}

\section{Applications to toric geometry}

\subsection{Toric stacks and Kawamata's Theorem.}

The applications in this section come from  toric geometry.
We will assume $k = \mathbb{C}$. 
Let $X, Y$ be quasi-smooth (i.e. simplicial) toric 
varieties with the action of the same split torus $T$, with two effective 
$T$-invariant $\mathbb{Q}$-divisors $B, C$ respectively, 
and assume that the coefficients of these divisors are 
in the set $\{1 - \frac{1}{r}| r \in \mathbb{N}\}$. To this data one can
associate two smooth Deligne-Mumford stacks $\ca X, \ca Y$, respectively, 
as in \cite{Kawamata:2005}. One reason why considering $B, C$ is useful is
that a quotient of a smooth variety $V$ by an effective 
action of a finite group $\Gamma$ may be smooth in codimension one yet
have non-trivial stabilizers there (automatically cyclic).  Then the multiplicities
of $B$ and $C$ will encode the information about the sizes of
such stabilizers.

If $B$ and $C$ are zero, the stacks $\ca X$, $\ca Y$ may be interpreted
in terms of the Batyrev-Cox quotient construction which we briefly 
recall.  
Let $\mathbb{P}(\Sigma)$ be a quasi-smooth  toric  variety
associated to a  simplicial fan $\Sigma$ which has $n$ one-dimensional
cones $\rho_1, \ldots, \rho_n$. 
Then, cf. e.g. \cite{CLS},  $\mathbb{P}(\Sigma)$
it can be realized as a quotient $(\C^n \setminus \mathbb{B})/G$ where 
the basis $u_1, \ldots, u_n$ in $\C^n$ is in bijective correspondence
with $\rho_1, \ldots, \rho_n$, $\mathbb{B}$ is a 
union of some coordinate subspaces of codimension at least 2, and $G$ is 
an algebraic subgroup of $(\C^*)^n$ with its natural action on $\C^n$.
Thus $G$ itself is isomorphic to a product of several copies of $(\C^*)$ 
and a finite abelian group. The assumption that $\Sigma$ is simplicial 
ensures, cf. \textit{loc. cit.}, that $G$-action on $\C^n \setminus \mathbb{B}$
has finite stabilizers and thus there exists a smooth Deligne-Mumford
quotient stack $\mathcal{P}(\Sigma)$ with coarse moduli space isomorphic to
 $\mathbb{P}(\Sigma)$. 

We will use a version of the Kawamata's theorem on equivalence of 
derived categories for toric stacks. In our setup, we assume that 
in the diagram \eqref{main} all stacks are toric with the action
of the same torus $T$, and that all morphisms are $T$-equivariant.

\begin{theorem}
\label{Ka}
In the situation described, there exists a $T$-equivariant kernel
$\mathcal{F}$ on $\ca X \times \ca Y$ given by a
direct image of a $T$-equivariant object in $D^b(\ca X \times_{Z} \ca Y)$,
such that the corresponding Fourier-Mukai transform
$\Phi_{\scr F}: D^b(\ca X) \to D^b(\ca Y)$ is an equivalence. 
\end{theorem}
\noindent
\textit{Proof.} By the relative version of the toric MMP
explained in \cite{FS} and a standard argument modeled
on the proof of Theorem 12.1.8 in \cite{Ma}, the 
birational isomorphism $X \to Y$ may be decomposed into a finite
sequence of $T$-equivariant divisorial contrations and flips \textit{over Z} 
which are log crepant. For every 
such contraction or flip we can apply Theorem 4.2 in \cite{Kawamata:2005}.
Observe that the kernels in \textit{loc. cit.} are indeed $T$-equiariant 
and given by direct 
images from $\ca X\times_{Z} \ca Y$. The convolution of several
such kernels (giving the composition of equivalences) also satisfies this
condition.
$\square$

\subsection{Hypersurfaces in simplicial toric varieties.}

As another application of Theorem 1.1 we extend a result of Orlov 
on the derived category of a Calabi-Yau hypersurface in a weighted 
projective space. From now on the log divisors $B, C$ are zero.
Our application is based on the following result due to M.U. Isik, 
cf. \cite{Isik:2009}. For completeness we reproduce a sketch of the proof
with the kind permission of the author. 

\begin{proposition}
\label{isik}
Let $X$ be a smooth Deligne-Mumford stack over $\C$ and
$s$ a regular section of a vector bundle $E$ with determinant dual to the 
canonical bundle of $X$. Let $Y \subset X$ be the zero scheme of $s$ 
and $Z \subset E^\vee$ the zero scheme of 
$s$ viewed as a fiberwise linear function
on the total space of the dual bundle $E^\vee$. Then $D^b(Y)$ is equivalent
to the Karoubian completion (or split completion) 
of the equivariant singular category $D_{sg}^{\C^*} (Z)$, which obtained by adding images of all projectors.
The action 
of $\C^*$ on $Z$ is restricted from the natural action by scalar dilations
on the fibers of $E^\vee$. 
\end{proposition}
\noindent
\textit{Proof.} The proof involves sheaves of graded DG-algebras (and
modules over them) which have homological degree $\deg_h$ and internal
degree $\deg_i$. First one identifies (up to derived equivalence) 
complexes of
graded $\mathcal{O}_Z$-modules with 
complexes of graded modules over the Koszul 
DG algebra 
$$
\mathcal{B} = (\varepsilon Sym^\bullet(E) \to Sym^\bullet(E)); \qquad
d(\varepsilon) = s
$$
where the formal variable $\varepsilon$ satisfies $\varepsilon^2 = 0$
and the differential is linear with respect to $Sym^\bullet(E)$. 
In both categories the differential increases $\deg_h$ by one and
preserves $\deg_i$. Note
that $\deg_h(\varepsilon)= -1$, $\deg_i (\varepsilon) = 1$, 
$\deg_h(Sym^\bullet (E)) = 0$ while the internal degree on $Sym^\bullet(E)$
is given by the usual polynomial degree.

Next use fiber-by-fiber BGG correspondence (or Koszul duality) to 
construct a derived equivalence between graded DG modules over $\mathcal{B}$
and graded DG-modules over the ``dual" algebra
$$
\mathcal{A} = \Lambda^\bullet(E^\vee) \otimes_{\mathcal{O}_X} 
\mathcal{O}_X[t]
$$
where the differential is $t$-linear and satisfies $d(f) = t \cdot \langle s, f
\rangle$ for a local section $f$ of $E^\vee$ (this admits a unique
extension to the exterior algebra by Leibniz rule). This time
$\deg_i(t) = \deg_i(f) = 1$, $\deg_h(t)= 2, \deg_h(E^\vee)=1$. 

The composition of these two equivalences sends perfect complexes of 
graded $\mathcal{O}_Z$-modules map to complexes of graded $\mathcal{A}$-modules
on which $t$ acts nilpotently. Therefore, using a version of Thomason's 
localization theorem \cite{TT}, Lemma 5.5.1, we conclude that the
split completion of $D^{\C^*}_{sg}(Z)$
is equivalent to the category of graded modules over the localization 
$\mathcal{A}_t$. Since $t$ is now invertible, multiplication by it
identifies all homogeneous components (with respect to $\deg_i$), 
i.e. the category is equivalent to 
non-graded modules over the Koszul resolution $\mathcal{A}'$ of $\mathcal{O}_Y$.
Using the derived equivalence between $\mathcal{A}'$-modules and
$\mathcal{O}_Y$-modules we obtain the result.
$\square$

\bigskip
\noindent
Now consider a projective simplicial toric stack as in Section 3.1. 
We  fix the fan 
$\Sigma$ and drop it from notation, writing simply $\ca P$ and 
$\mathbb{P}$. We will require  that $K_{\ca P}^\vee$ is nef, i.e. a positive
power of the anti-canonical bundle on $\ca P$ descends to 
a nef bundle on $\mathbb{P}$.

The Picard group of the stack $\mathcal{P} = [(\C^n \setminus \mathbb{B})/G]$ 
can be identified with the
$A = Hom(G, \C^*)$ since $\mathbb{B}$ has codimension $\geq 2$ and therefore
any line bundle on $\C^n \setminus \mathbb{B}$ is trivial. 
Explicitly, for $a \in A$ we have a line bundle $L_a$ over $\CP$ with 
the total space $[(\C^n \setminus \mathbb{B}) \times \C]/G$ where
 $G$ acts on the
second factor by the character $-a$. The sections of $L_a$, viewed
as regular functions on the total space of $L_{-a}$ 
linear along the fibers, may be identified with the space 
$\C[x_1, \ldots, x_n]^a$ of polynomials on $\C^n$, on which $G$ acts via
the character $a$. Here $x_1, \ldots, x_n$ is the basis dual 
to the basis $u_1, \ldots, u_n$ of Section 3.1. Observe that
$\C[x_1, \ldots, x_n]^a$ is non-zero if and only if $a$ belongs to
the semigroup $A^+ \subset A$ formed by all non-negative integral
linear combinations of the $G$-weights $a_1, \ldots, a_n$ of 
$x_1, \ldots, x_n$, respectively. Note that all elements in $A^+ \setminus 0$
have infinite order if $\mathbb{P}$ is projective, otherwise 
some monomial with non-negative exponents would give a non-constant
regular function on $\mathbb{P}$.
The line bundle $L_a$
descends to the coarse moduli space $\mathbb{P}$ precisely when 
$a$ is trivial on all stabilizers of the $G$-action on $\C^n \setminus B$, 
hence the Picard group of $\mathbb{P}$ is a subgroup of finite index in $A$. 
We also recall that
$K_{\mathcal{P}(\Sigma)}^\vee$ is isomorphic to $L_{c}$ with
$$
c = a_1+ \ldots + a_n.
$$
In the smooth case this is proved in Section 4.3 of \cite{Fu}
which is sufficient for us since $\mathbb{P}$ is smooth in codimension 1.
Choose and fix a polynomial $f \in \C[x_1, \ldots, x_n]^c$ and consider
its zero locus $\ca M$, a closed substack of $\CP$ which has trivial 
canonical bundle by the adjunction formula. Note that $\ca M$ may be
identified with the quotient stack $\big[Z(f) \cap (\C^n \setminus Z)\big]/G$.

Let $H \subset G$ be the kernel of the character 
$c: G\to \C^*$. Then the group of characters of $H$ may be identified with 
$A_H= A/\mathbb{Z}\cdot c$ and we denote by $A^+_H$ the image of $A^+$ in 
$A_H$. The $G$-action on $\C^n$ restricts to $H$ and for a generic
choice of $a \in A^+_H$ the $H$-linearization of the trivial bundle, induced
by $a$, will satisfy the property that all semi-stable points are 
stable (in fact, it suffices to require that $a$ does not belong to 
any subgroup in $A_H$ which is generated by a finite subset in 
the image of $\{a_1,\ldots, a_n\}$ and has non-maximal rank). 
Let $\ca Y^a = (\C^n)^s/H$ be the corresponding 
quotient stack; its coarse moduli space $Y$ is just the GIT quotient with
respect to the linearization induced by $a$. 

Since the character $c: G \to \C^*$ is trivial on $H$, the polynomial
$f$ on $\C^n$ is $H$-invariant  and therefore
descends to a morphism $g: \ca Y^a \to \mathbb{A}^1$. For the
same reason, i.e. triviality of $c$ on $H$, 
the stack $\ca Y^a$ has trivial canonical bundle. 
The action of $(\C^*)^n$ on $\C^n$ descends to the action of
the torus $T_H = (\C^*)^n/H$ on 
$\ca Y^a$. Observe that there exists a short exact sequence
$$
1 \to \C^* \to T_H \to T \to 1
$$ 
where the subgroup $\C^*$ may be identified with $G/H$.

\begin{corollary}  In the above setting, the derived category 
$D^b(\ca M)$ is equivalent to the split
completion of the
equivariant singular category $D^{\C^*}_{sg} (g^{-1}(0))$. 
\end{corollary}
\begin{proof} First, we can view $f$ as a function, which we denote by the 
same letter, on the total space of the canonical bundle $K_{\CP}$.
By Proposition \ref{isik} the derived category of $\ca M$ is equivalent
to the split completion of the equivariant singular category 
 of the zero fiber of $f: K_{\CP} \to \mathbb{A}^1$. Observe that
both $K_{\CP}$ and $\mathcal{Y}^a$ are toric with the same torus 
$T_H$. For $K_\CP$ this follows from the fact that it is 
a quotient of an open subset in $\C^n \times \C$ by $G$, and if
we use the character $c$ to lift the natural embedding 
$G\subset (\C^*)^n$ to an embedding $G\subset (\C^*)^n \times \C$, then
the quotient $(\C^*)^n \times \C/G$ is canonically isomorphic to $T_H$. 

Since both $K_\CP$ and $\ca Y^a$ have trivial canonical bundles, by 
Theorem 1.1 and Theorem \ref{Ka} it suffices to find 
an affine variety $Z$ which is toric with respect to  $T_H$,
and two proper equivariant birational morphisms $\phi: K_\CP \to Z$, 
$\psi: \ca Y^a \to Z$ which induce projective morphisms on coarse moduli
spaces. Since generic stabilizers of $K_\CP$ and $\ca Y^a$
are trivial and the morphisms to the coarse moduli spaces 
are proper (\cite[Prop. 2.11]{Vi}, \cite[Prop. 4.2]{Ed}), it suffices to construct birational projective
toric morphisms from the coarse moduli spaces to $Z$.

Choose $Z$ to be the spectrum of the ring of invariants
$\C[x_1, \ldots, x_n]^H$. Then the coarse moduli space $Y^a$ of $\ca Y^a$
is projective over $Z$: by GIT it is a $Proj$ of a graded ring
$\bigoplus_{l \geq 0} R_l$ with $Z \simeq Spec(R_0)$. 
On the other hand, the moduli space of $K_\CP$ is the geometric
quotient of an open subset in $\C^n\times \C$ by $G$, where the action on
the second factor is via the character $c$ of $G$. 

Consider the ring of invariants $\C[x_1, \ldots, x_n, z]^G$ for this 
action. Since $z$ is dual to the last coordinate vector in 
$\C^n \times \C$, the character of the $G$-action on it is  $(-c)$. 
Every $G$-invariant polynomial is a linear combination of terms
of the form $h(x) z^l, l \geq 0$ where $G$ acts on 
$h(x) \in \C[x_1, \ldots, x_n]$ via the character $l\cdot c$. 
In particular, each $h(x)$ is $H$-invariant. Therefore evaluation
$z \mapsto 1$ gives a ring homomorphism
$$
\C[x_1, \ldots, x_n, z]^G\to \C[x_1, \ldots, x_n]^H
$$
To show that it is as isomorphism observe that  every $H$-invariant 
polynomial is a sum of $H$-invariant monomials and that of each 
$H$-invariant monomial $x^\alpha$ the group $G$ acts via a character
with is trivial on $H$.  Since 
$Ker (A \to A_H) = \mathbb{Z}\cdot c$, such a character is 
a multiple  of $c$. But a multiple $l\cdot c$ with  \textit{negative} $l$ 
may not be a $G$-weight of 
monomial $x^\alpha$ with \textit{positive} exponents, otherwise 
$x^\alpha (x_1 \ldots x_n)^{-l}$ would give a non-constant regular
function on the projective variety $\mathbb{P} = (\C^n \setminus \mathbb{B})/G$. 
We proved that 
$$
\C[x_1, \ldots, x_n]^H = \bigoplus_{l \geq 0}\C[x_1, \ldots, x_n]^{l \cdot c}
$$
which implies that the above map of invariants is an isomorphism.
We also observe that the grading on the left hand side is precisely the
grading resulting from the action of $\C^*\subset T_H$ on 
$Z$. 

It remains to show that the coarse moduli space of $K_\CP$, i.e. the
geometric quotient $\big((\C^n \setminus \mathbb{B}) \times \C\big)/G$, 
admits a projective
morphism to $Z$. Since the stabilizers of the $G$-action are
finite the assertion would follow from GIT
if  $(\C^n \setminus \mathbb{B}) \times \C$ 
is the set of semistable points for
some $G$-linearization of the trivial bundle, coming from a character $a \in A$.
Choose any $a$ which gives an ample line bundle in $Pic(\mathbb{P})
\subset A$. 
We need to prove that $(\C^n \setminus \mathbb{B}) \times \C$
is precisely the set of points for which one can find a non-vanishing
polynomial $f(x) z^l$ of $G$-weight $m \cdot a$, $m \geq 0$. 

For any point $(p, q) \in (\C^n\setminus B)\times \C$, $p$ will project
to a point $\overline{p} \in \mathbb{P}$ and by ampleness there is 
a section of $L_a^{\otimes m}$ not vanishing at $\overline{p}$.
This section gives a polynomial $f(x)$ of $G$-weight $a \cdot m$
non-vanishing at $p$, hence $f(x) z^0$ is a quasi-invariant polynomial 
not vanishing at $(p, q)$. 

On the other hand, we want to show that any $f(x) z^l$ of $G$-weight
$m \cdot c$ will vanish at any point $(p, q) \in \mathbb{B} \times \C$.
It suffices to show that any polynomial $f(x)$ of weight 
$a' = l \cdot c + m \cdot a$ with $m > 0, l \geq 0$, will vanish at $p \in B$. 
Since we assumed that $K^\vee_\CP$ is nef, replacing
$f(x)$ by its positive power  we can assume that $a'$ is an ample class
in $Pic(\mathbb{P})$. We can also assume that
$f(x)$ is a monomial $x^\alpha$. 
 
Now take a closer look at $\mathbb{B}$, a union of coordinate subspaces 
of the form
$\mathbb{B}_j = \{x_i = 0| i \in P_j\}$ where $P_j \subset \{1, \ldots, n\}$
is a subset called a \textit{primitive collection} (and the index $j$ will
run over all primitive collections). See Section 5.1 
of \cite{CLS} for details. 
The crucial observation is that each $P_j$ gives a class $r_j$
in the cone of effective curves on $\mathbb{P}$. Since $a'$ is ample
and the cone of effective curves is spanned by $r_j$, see 
Theorem 6.3.10 in \textit{loc. cit}, 
the intersection number $a' \cdot r_j$ should be positive for all $j$, which will imply
that for any $j$ there exists $i \in P_j$ which gives a positive 
exponent $\alpha_i$ in the monomial $x^\alpha = x_1^{\alpha_1} \ldots x_n^{\alpha_n}$. Then $x^\alpha$ must vanish on each $B_j$, which will
imply the assertion about the set of stable points. 

In more detail: each variable $x_i, i = 1, \ldots, n$ corresponds
to a  one-dimensional cone 
$\rho_i \in N_{\mathbb{R}}$ in the
fan $\Sigma$ defining our toric variety, and hence to a torus-invariant
prime divisor $D_i$ in $\mathbb{P}$. One the other hand, the space of
numerical classes of curves in $\mathbb{P}$ can be identified with 
the space of vectors $C= (\beta_1, \ldots, \beta_n)$ such that
$\sum b_i \rho_i = 0$ in $N_{\mathbb{R}}$, in such a way that
for a divisor $D= \sum \alpha_i D_i$ the intersection number $D\cdot C$
is given by $\sum \alpha_i \beta_i$, cf. Exercise 6.3.3 in \textit{loc. cit.}

Now, for any primitive collection $P_j$, Batyrev's construction gives
a relation in $N_{\mathbb{R}}$ of the type
$$
r_j = \sum_{i \in P_j} \rho_i - \sum_{s=1}^n c_s \rho_s 
$$
with $c_s \in \mathbb{Q}_{\geq 0}$, cf. Definition 6.3.9
of \textit{loc. cit.} Moreover, since $r_j$ corresponds to 
a numerically effective class of curves by ampleness of $a'$ we
should have 
$$
\sum_{i \in P_j} \alpha_i - \sum_{s = 1}^n c_s \alpha_s > 0.
$$
All $c_s$ and $\alpha_s$ are non-negative hence $\alpha_i$ 
for some $i \in P_j$. That is, the monomial 
$x^\alpha$ contains the variable $x_i$ with a positive exponent and hence
vanishes on the subspace $\mathbb{B}_j \subset \mathbb{B}$. 
Since this holds for every
primitive collection $P_j$, and $\mathbb{B}$ is the union of 
$\mathbb{B}_j$, we conclude
that $x^\alpha$ vanishes on $\mathbb{B}$. 

\bigskip
\noindent
To sum up: we have proved that 
$\ca X = K_\CP$ and  $\ca Y^a$ admit
projective toric morphisms (automatically birational) onto the affine
variety $Z= Spec(\C[x_1, \ldots, x_n]^H)$. Since the functions on 
$\ca X$ and $\ca Y^a$ are pulled back from $Z$, we can apply 
Theorem \ref{Ka} and
Theorem 1.1 to conclude that the equivariant singular categories 
of $\ca X_0$ and $\ca Y^a_0$ are equivalent, which finishes the proof. 
\end{proof}
 
\begin{example}
\cite[Thm. 3.12, Calabi-Yau case]{Orlov:2005}  Take the weighted projective space $\mathbb{P} = \mathbb{P}(\overline{a}):=\mathbb{P}(a_0,\cdots, a_n)$ with $a_i>0$ for all $i$.  Let $f$ be a quasi-homogeneous polynomial that is invariant under the action of $H = \Z_N$ where $N=\sum a_i$ and $\ca M$ the zero set of $f$
in the toric stack corresponding to $\mathbb{P}(a)$.  
Denote by $g$ the induced function on the quotient stack $\ca Y = [\C^n/\Z_N]$
(here $H$ is finite, and no choice of $a$
is needed). Then there is an equivalence between $D^b(\ca M)$ and the split completion of $D^{\C^*}_{sg}(g^{-1}(0))$.
\end{example}

\begin{example} Take $\ca P = \mathbb{P} = \mathbb{P}^1\times \mathbb{P}^1$ 
with $A = \Z \oplus \Z$. The weights of the homogeneous 
coordinates $x_1, \ldots,
x_4$ are $a_1 = a_2 = (1, 0)$ and $a_3 = a_4 = (0, 1)$. A quick calculation
shows that $H \simeq \C^* \times \Z_2$ hence $\ca Y^a$ is a GIT quotient
of $\C^4$ by $H$. Restricting the action on $\C^4$ to $\C^* \subset H$ we find 
that the weights are $(+1, +1, -1, -1)$. There are two essentially different
linearizations: one which gives a positive weight when restricted to 
$\C^*$ and one which gives a negative weight.
The two resulting quotients $\ca Y^+$ and $\ca Y^-$ differ by a 
standard toric flop. In particular, they have equivalent derived categories 
and applying Theorem 1.1 with $\ca X = \ca Y^+$, $\ca Y= \ca Y^-$
we see that the equivariant singular derived categories 
of $(g^+)^{-1}(0)$ and $(g^-)^{-1}(0)$ are equivalent. Their split
completions are further equivalent to the derived category of the 
elliptic curve $\ca M$ in $\mathbb{P}^1 \times \mathbb{P}^1$ by 
Corollary 3.3.
\end{example}

\begin{example}  Take $\ca P = \mathbb{P} =\mathbb{F}_1$, 
the Hirzebruch surface, which 
may be identified with the quotient
$(k^4\setminus \mathbb{B})/(\C^*)^2$ where $\mathbb{B}$ is the union of two 
coordinate planes $x_1 = x_2 = 0$ and $x_3 = x_4 = 0$. The weights of 
the action of $G= (\C^*)^2$ on $\C^4$ are $(1, 0), (1, 1), (0, 1), (0, 1)$.
The subgroup $H\simeq \C^*$ is given by $(t^3, t^{-2})$ and hence the 
weights of the $H$ action on $\C^4$ are $(3,1, -2, -2)$. Again we have two 
linearizations for the $H$ action and the corresponding choices of
$\ca Y^a$  differ by a toric flop.  
 \end{example}    

\subsection{Products of hypersurfaces}

Assume that for $i = 1, \ldots m$, we are given a Calabi-Yau
hypersurface $\ca M_i$ in a projective simplicial toric stack 
$\ca P_i$ with nef anticanonical bundle. By adjunction, $\ca M_i$
is defined by vanishing of a section $t_i$ of the
anticanonical bundle $K_{\ca P_i}^\vee$. 
 We can view the product $\ca M = \ca M_1 \times 
\ldots \ca M_m$ as a complete intersection in 
the stack $\CP= \CP_1 \times \ldots \times \CP_m$. 

On the other hand, let $\ca Y^{a_i}$ be a stack corresponding to 
a character $a_i$ as in Section 3.2, and $g_i: \ca Y^{a_i} \to \mathbb{A}^1$
the function induced by the defining equation of $\ca M_i$. This gives
a function $g = g_1 + \ldots + g_m$ on the direct product $\ca Y = \ca Y^{a_1}
\times \ldots \times \ca Y^{a_m}$. The action of $\C^*$ on each 
$\ca Y^{a_i}$ induces the diagonal action on $\ca Y$.
\begin{corollary}
There exists an equivalence between $D^b(\ca M)$ and the split 
completion of $D^{\C^*}_{sg} (g^{-1}(0))$.
\end{corollary} 
\begin{proof}
Let $\ca X= K_{\ca P_1} \times \ldots K_{\ca P_m}$ be the product
of total spaces of canonical bundles on the stacks, and 
$f: \ca X\to \mathbb{A}^1$ the fiberwise linear function corresponding 
to $t_1 \oplus \ldots \oplus t_m$. By the same argument 
as with a single hypersurface, both $\ca X$ and $\ca Y$ admit
a proper birational morphism to the product affine variety 
$Z = Z_1 \times \ldots \times Z_m$ which induces a projective morphism on 
moduli spaces. Thus $D^{\C^*}_{sg}(f^{-1}(0)) \simeq D^{\C^*}_{sg} (g^{-1}(0))$
as before. Now we apply Proposition \ref{isik} to the direct sum of
anticanonical bundles on $\ca P_1 \times \ldots \times \ca P_m$.
\end{proof} 
\begin{remark}
One could attempt to apply the same reasoning to general Calabi-Yau
complete intersections in a toric variety $\ca P$ but it leads to a 
technical difficulty. Let $\beta_1, \ldots, \beta_m$ be elements in 
the stack Picard group $A$, such that
$$
\beta_1 + \ldots + \beta_m = a_1 + \ldots + a_n
$$
Assume in addition that each bundle $L_{\beta_i}$ has a non-zero 
regular section $t_i$ and that $t_1, \ldots, t_m$ form a regular sequence.
This implies, in particular, that each $\beta_j$ is a linear combination 
of the $a_i$ with non-negative coefficients. Let $A_H$ be the 
quotient by the subgroup of $A$ spanned by $\beta_1, \ldots, \beta_m$
and $H \subset G$ the dual group. 
Suppose we wanted to establish a derived equivalence between the total
space $\ca X$ of the bundle $L_{-\beta_1} \oplus \ldots \oplus L_{-\beta_m}$ 
and a quotient $\ca Y$ of some open subset of $H$-stable points 
in $\C^n$, by the action of $H$. In particular, we would want  
proper birational morphisms from $\ca X$ and $\ca Y$ to the same 
variety (or stack) $Z$.

If we define $Z$ as a spectrum of a ring of invariants, 
then for $\ca Y$ this ring should be $\C[x_1, \ldots, x_n]^H$,
while for $\ca X$ the ring is $\C[x_1, \ldots, x_n, z_1, \ldots, z_m]^G$ 
where $G$ acts on the extra variable $z_j$ via the character $-\beta_j$.
Sending each $z_j$ to $1$ (or to some nonzero constant $c_j$) we will get
a map 
$$
\C[x_1, \ldots, x_n, z_1, \ldots, z_m]^G \to \C[x_1, \ldots, x_m]^H.
$$
We observe that this map is injective if and only if $\beta_1, \ldots, \beta_m$
are linearly independent over $\Z$, which is a reasonable condition. However, 
the surjectivity is less trivial: one would need to ensure that for
every equality in $A$
$$
p_1 a_1 + \ldots + p_n a_n = q_1 \beta_1 + \ldots q_m \beta_m
$$
the condition $p_i \geq 0$ for all $i$ implies the condition  $q_j \geq 0$
for all $j$. The reader is invited to check that this fails
e.g. when $\ca P = \mathbb{P}^2 \times \mathbb{P}^2$ with $Pic = \Z\oplus \Z$, 
and 
$$
a_1 = a_2 = a_3 = (1, 0); \quad a_4 = a_5 = a_6 = (0,1); \quad
\beta_1 = a_1 + a_4 + a_5; \quad \beta_2 = a_2 + a_3 + a_6.
$$
In such a situation 
$\ca X$ and $\ca Y$ would be proper over non-isomorphic, although birational, 
affine varieties.
Perhaps one could modify our construction, e.g. by passing to a partial 
compactification of $\ca X$ or $\ca Y$, to obtain the analogue
of Corollary 3.7 in  a more general situation.
\end{remark}

\end{document}